\documentclass[reqno]{amsart}
\usepackage{amsmath, amsthm, amssymb, amscd, hyperref}
\usepackage{mathabx} 
\usepackage{xypic}
\xyoption{all}
\usepackage{enumitem}
\setenumerate{label=\textnormal{(\arabic*)}} 
\usepackage{harpoon} 

\usepackage{mathrsfs} 
\usepackage{wasysym} 

\usepackage{color}   

\newcommand*\colvec[3][]{
    \begin{pmatrix}\ifx\relax#1\relax\else#1\\\fi#2\\#3\end{pmatrix}
} 

\makeatletter 
\let\@wraptoccontribs\wraptoccontribs 
\makeatother

\theoremstyle{plain}
\newtheorem*{KStheorem}{Kochen-Specker Theorem}
\newtheorem{theorem}{Theorem}
\newtheorem{corollary}[theorem]{Corollary}
\newtheorem{lemma}[theorem]{Lemma}
\newtheorem{proposition}[theorem]{Proposition}
\newtheorem{question}[theorem]{Question}

\theoremstyle{definition}
\newtheorem*{definition*}{Definition}

\DeclareMathOperator{\pSpec}{\mathit{p}-Spec}

\DeclareMathOperator{\Proj}{Proj}
\DeclareMathOperator{\Idpt}{Idem}

\DeclareMathOperator{\diag}{diag}

\DeclareMathOperator{\rad}{rad}
 
\newcommand{\C}{\mathbb{C}}
\newcommand{\Q}{\mathbb{Q}}
\newcommand{\Z}{\mathbb{Z}}
\newcommand{\F}{\mathbb{F}}
\newcommand{\R}{\mathbb{R}}
\newcommand{\M}{\mathbb{M}}

\newcommand{\setS}{\mathcal{S}}

\renewcommand{\O}{\mathcal{O}}

\newcommand{\comm}{\odot}
\newcommand{\sa}{\mathrm{sa}}
\newcommand{\sym}{\mathrm{sym}}

\newcommand{\separate}{\bigskip}

\begin{document}

\title[Integer Kochen-Specker contextuality]{Minimal ring extensions of the integers exhibiting Kochen-Specker contextuality}

\author{Ida Cortez}
\address{Bowdoin College\\
Department of Mathematics\\
8600 College Station\\
Brunswick, ME 04011}
\email{idac@live.com}

\author{Camilo Morales}
\address{Harvey Mudd College\\
Department of Mathematics\\
320 E.~Foothill Blvd.\\
Claremont, CA 91711, USA}
\email{camorales@g.hmc.edu}

\author{Manuel Reyes}
\address{University of California, Irvine\\
Department of Mathematics\\
340 Rowland Hall,\\
Irvine, CA 92697-3875, USA}
\email{manny.reyes@uci.edu}
\urladdr{https://math.uci.edu/~mreyes/}
\thanks{Cortez was supported by a Maine Space Grant Consortium Fellowship. Morales was supported by the Summer Undergraduate Research Fellowship program at UC Irvine. Reyes was supported in part by National Science Foundation grant DMS-2201273.}

\date{November~19, 2025}

\begin{abstract}
This paper is a contribution to the algebraic study of contextuality in quantum theory. 
As an algebraic analogue of Kochen and Specker's no-hidden-variables result, 
we investigate rational subrings 
over which the partial ring of $d \times d$ symmetric matrices ($d \geq 3$) admits no morphism to a nonzero commutative ring, which we view as an ``algebraic hidden state.'' 
For $d = 3$, the minimal such ring is shown to be $\Z[1/6]$, while for $d \geq 6$ the minimal subring is $\Z$ itself.
The proofs rely on the construction of new sets of integer vectors in dimensions~3 and~6 that have no Kochen-Specker coloring.
\end{abstract}

\dedicatory{In fond memory of Professor S.~Tariq Rizvi.}

\maketitle

\section{Introduction}
\label{sec:intro}

Quantum contextuality~\cite{BCGKL} is a foundational feature of quantum theory, which has grown to a widely studied phenomenon since the work of Bell~\cite{Bell} and Kochen-Specker~\cite{KochenSpecker}.  
Recall that the observables of a quantum-mechanical system are represented by elements of a noncommutative algebra of operators. 
For a pair of observables that do not commute, the uncertainty principle forbids us from simultaneously measuring the values of these observables with perfect certainty. Conversely, a pair of observables is called \emph{commeasurable} if they commute, so that there is no such restriction on their simultaneous measurement. Broadly speaking, \emph{contextuality} in quantum mechanics is the principle that the value of a measured observable is dependent on the commeasurable observables that happen to be probed by a measurement apparatus.

Contextuality has been studied through a number of different mathematical lenses. In recent years this has included sheaf theory~\cite{AbramskyBrandenburger, AMB}, Boolean algebra~\cite{BMR, AbramskyBarbosa}, and topology~\cite{ORBR, OkayRaussendorf, OkaySheinbaum}. 
In this paper we take inspiration directly from Kochen and Specker's initial treatment of contextuality in terms of partial algebras, leading to a ring-theoretic and number-theoretic investigation as we explain below.

\separate

In the foundations of quantum physics, a \emph{hidden variable theory} is an attempt to reduce the uncertainty of the outcomes of quantum measurements to a statistical average over classical ``hidden variables'' that behave deterministically.
Kochen and Specker's investigation of contextuality centered on the proof that certain noncontextual hidden variable theories are impossible to construct. Their formulation of noncontextual hidden variable theories has a strongly algebraic flavor. 
We recall some terminology from~\cite{KochenSpecker} in order to state their result.

Let $K$ be a commutative ring. A \emph{partial $K$-algebra} is a set $A$ equipped with a reflexive, symmetric binary relation $\comm \subseteq A \times A$ called \emph{commeasurability}, along with partially defined operations $+,\, \cdot \ \colon \odot \to A$ of addition and multiplication as well as a globally defined scalar multiplication $K \times A \to A$, and distinguished elements $0,1 \in A$, subject to the following axioms:
\begin{enumerate}[label=(P\arabic*)]
\item $0$ and $1$ are commeasurable with all elements of $A$;
\item the partial binary operations are commutative when defined: if $a, b \in A$ with $a \odot b$ then $a+b = b+a$ and $ab = ba$;
\item all operations preserve commeasurability: if $a,b,c \in A$ are pairwise commeasurable, then $(a+b) \odot c$, $ab \odot c$, and $\lambda a \odot c$ for all $\lambda \in K$;
\item if $a,b,c \in A$ are pairwise commeasurable, then the set of evaluations of all polynomials $\{f(a,b,c) \in A \mid f(x,y,z) \in K[x,y,z]\}$ forms a commutative $K$-algebra under the restricted operations from $A$.
\end{enumerate}
Assuming (P1)--(P3), condition~(P4) is equivalent to the statement that every set $S \subseteq A$ of pairwise commeasurable elements is contained in a set $C \subseteq A$ of pairwise commeasurable elements such that the restriction of all operations makes $C$ into a commutative $K$-algebra with additive and multiplicative identities $0$ and $1$.

A \emph{morphism} of partial $K$-algebras $f \colon A \to B$ is a function that satisfies the following conditions for any $a, b \in A$:
\begin{enumerate}[label=(M\arabic*)]
\item $f(0) = 0$ and $f(1) = 1$;
\item $f(\lambda a) = \lambda f(a)$ for all $\lambda \in K$;
\item $a \odot b \implies f(a) \odot f(b)$, $f(a+b) = f(a) + f(b)$, and $f(ab) = f(a) f(b)$.
\end{enumerate}
Evidently, if $C \subseteq A$ is a commeasurable subalgebra, then $f$ (co)restricts to a homomorphism $C \to f(C)$ of commutative $K$-algebras.

\separate

Kochen and Specker noted that the set of observables $\O$ of a given physical system forms a partial $\R$-algebra. For instance, a $d$-level quantum system has observables $\O = \M_d(\C)_\sa$ consisting of self-adjoint $d \times d$ complex matrices. If we define commeasurability to be commutativity of operators as described above, then $\O$ forms a partial $\R$-algebra under addition and multiplication of commuting pairs of operators, and the subset $\M_d(\R)_\sym$ of real symmetric matrices forms a partial subalgebra. More generally, the set of self-adjoint elements of any complex $*$-algebra of operators on a Hilbert space evidently forms a partial $\R$-algebra, with commeasurability again defined by commutativity.

A \emph{Kochen-Specker hidden variable theory} for a system can now be defined concisely as an injective morphism of partial $\R$-algebras
\[
h \colon \O \to \R^\Omega,
\]
where $\O$ is the partial $\R$-algebra of observables on the system, $\Omega$ is a set (viewed as a space of \emph{hidden states}~\cite[p.~62]{KochenSpecker} on which the hidden variables are defined), and $\R^\Omega$ is the algebra of all real-valued functions on $\Omega$. 
Note that evaluation at each element $x \in \Omega$ composes with $h$ to produce a morphism of partial $\R$-algebras $\psi_x \colon \O \to \R$, which assigns a definite value to every observable in a way that preserves algebraic relations between commeasurable observables. Conversely, if sufficiently many such \emph{hidden pure state} morphisms $\psi \colon \O \to \R$ of partial algebras can be produced to separate distinct observables, we may take the product over these morphisms to construct a hidden variable theory $h \colon \O \to \R^\Omega$.

Kochen and Specker proved~\cite[pp.~66 \&~70]{KochenSpecker} the following obstruction to any such hidden variable theory.

\begin{KStheorem}[\cite{KochenSpecker}]
Let $d \geq 3$ be an integer. There is no morphism of partial $\R$-algebras from $\M_d(\R)_\sym$ to any nonzero commutative $\R$-algebra. Thus there is no Kochen-Specker hidden variable theory for $\O = \M_d(\C)_\sa$.
\end{KStheorem}

This was achieved by showing the non-existence of certain colorings on the partial Boolean algebra of projections, which in turn was deduced from the uncolorability of vectors that represent rank-1 projections. (The definition of these colorings and of partial Boolean algebras will be recalled in later sections of this paper.) 
Furthermore, they exhibited~\cite[Section~6]{KochenSpecker} a hidden variable theory in the case of dimension~2 by constructing an injective morphism of partial $\R$-algebras $\M_2(\C)_\sa \to \R^\Omega$ where $\Omega = S^2 \subseteq \R^3$ is the unit sphere.

\separate

The formulation of the Kochen-Specker Theorem suggests a purely algebraic treatment of contextuality in the setting of noncommutative algebra via partial rings, explored in~\cite{Reyes, BMR}. It was shown in~\cite{Reyes} (see also~\cite{Reyes:tour}) that Kochen and Specker's original result provides an obstruction to the existence of functors that extend the prime spectrum from commutative to noncommutative rings, in the sense that such functors assign the empty set to the rings $\M_d(\C)$ for $d \geq 3$. 
These two results were connected by the spectrum of \emph{partial prime ideals}. 

The problem of extending this result to matrix rings over any base naturally led to a more detailed study of partial algebras in~\cite{BMR}. By a \emph{partial ring} we simply mean a partial $\Z$-algebra. In that paper, the obstruction for spectrum functors was reduced to $\M_d(\Z)$. However, it was also shown that symmetric integer matrices have an algebraic analogue of a Kochen-Specker hidden state, taking values in a finite field, described as follows.

\begin{definition*}
Let $R$ be a partial ring. An \emph{algebraic hidden state (with values in $C$)} for $R$ is a morphism of partial rings $\psi \colon R \to C$ to a nonzero commutative ring $C$. 
\end{definition*}

For a commutative ring $K$, we let $\M_d(K)_\sym$ denote the set of symmetric $d \times d$ matrices over $K$; this forms a partial $K$-algebra with commeasurability again being commutativity, and all operations restricted from those of the full matrix algebra.
In the case $d = 3$, it was proved in~\cite{BMR} that:
\begin{itemize}
\item there exist morphisms of partial rings from $\M_3(\Z)_\sym$ to the finite fields $\F_{2^6}$ and $\F_{3^6}$;
\item there is no morphism of partial rings from $\M_3(\Z[1/30])_\sym$ to any nonzero commutative ring.
\end{itemize}
Intuitively speaking,
although the $3 \times 3$ symmetric matrices defined over $\Z$ have an algebraic hidden state, contextuality (interpreted as the absence of such ``hidden states'') emerges once we extend the base ring to $\Z[1/30]$.

This suggests the question of whether there are proper subrings of $\Z[1/30]$ for which  
algebraic hidden states are 
forbidden. In fact, the smallest such subring in dimension $d = 3$ is $\Z[1/6]$, according to the main result of this paper. 
Furthermore, for $d \geq 6$ contextuality is already present for the partial ring of symmetric integer matrices.
In dimensions~4 and~5, we do not have a final answer, but the our results suggest that the minimal contextual ring is $\Z[1/2]$.

\begin{theorem}\label{thm:map to commutative}
Let $N$ and $d$ denote positive integers.
\begin{enumerate}
\item The partial ring $\M_3(\Z[1/N])_\sym$ has no algebraic hidden states if and only if $N \equiv 0 \pmod{6}$.
\item If $d = 4,5$ and $N$ is even, then $\M_d(\Z[1/N])_\sym$ has no algebraic hidden states.
\item For $d \geq 6$, the partial ring $\M_d(\Z)_\sym$ has no algebraic hidden states.
\end{enumerate}
\end{theorem}

The proof of this result depends on the existence of uncolorable sets of integer vectors. For this reason, Section~\ref{sec:uncolorable} is devoted to the study of Kochen-Specker colorability of certain sets of integer vectors. 
The case of dimension $d = 3$ relies on a set constructed in~\cite{Bub:Schutte} along with a new set constructed in Theorem~\ref{thm:Q}. The case of $d \geq 4$ makes use of uncolorable sets in dimension~4 as in any of the references~\cite{CEG, Kernaghan, Peres:KS}.
The proof for $d \geq 6$ depends on a set whose uncolorability is verified by computational methods.

In Section~\ref{sec:rings} we apply the uncolorability of these vector sets to understand exactly which values of $N$ yield partial rings $\M_3(\Z[1/N])_\sym$ that have no prime partial ideals (Theorem~\ref{thm:prime}). 
This result on projection colorings is then applied to prove Theorem~\ref{thm:map to commutative}. 
We also show that certain partial Boolean algebras of projection matrices over various rings of number-theoretic interest are uncolorable (Corollary~\ref{cor:projections}).

\subsection*{Acknowledgments} 
This work would not have been possible without the previous efforts of Nikhil Dasgupta, Michael Pun, Samuel Swain, and Ira Hanzhao Li, and we are grateful for their contributions. We also thank the anonymous referee for some helpful corrections.

\section{Kochen-Specker colorings of integer vectors}
\label{sec:uncolorable}

Kochen and Specker famously reduced the proof of their no-hidden-variables result to the construction of a set of vectors in a Hilbert space that cannot be assigned a certain type of coloring, whose definition we now recall.
Let $S$ be a set of nonzero vectors in an $d$-dimensional Hilbert space $H \cong \C^d$. For the purposes of this paper, a \emph{Kochen-Specker (KS) coloring} of $S$ is $\{0,1\}$-coloring (i.e., a function $S \to \{0,1\}$) satisfying the following conditions:
\begin{enumerate}[label=(KS\arabic*)]
\item If $v, \lambda v \in S$ are collinear vectors (for $\lambda \in \C \setminus \{0\}$), then $v$ and $\lambda v$ are assigned the same color;
\item for any orthogonal set of vectors $v_1,\dots,v_m \in S$, at most one of the $v_i$ is colored~1;
\item if $m = d$ in (KS2), then exactly one of the $v_i$ is colored~1. 
\end{enumerate}
(This is slightly more general than some definitions in the literature, which require all vectors in $S$ to belong to an orthogonal basis in $S$.)
Given a subset $S_0 \subseteq S$, a Kochen-Specker coloring of $S$ restricts to a coloring of $S_0$; thus if $S_0$ has no Kochen-Specker coloring then also $S$ is uncolorable.

Kochen and Specker's no-hidden-variables theorem was proved by demonstrating the existence of a set of~117 three-dimensional vectors for which there is no KS coloring. Since then, a wide variety of KS~uncolorable vector sets of varying dimensions have been found. 
The proof of Theorem~\ref{thm:map to commutative} relies on the existence of Kochen-Specker uncolorable sets of \emph{integer} vectors, i.e.~vectors whose entries are all integers. 
Some of the best-known KS vector sets in fact consist of integer vectors.
This includes the collections of vectors in dimension~4 with entries in $\{-1, 0, 1\} \subseteq \Z$, given by Peres~\cite{Peres:KS}, Kernaghan~\cite{Kernaghan}, and Cabello, Estebaranz, and Garc\'{i}a-Alcaine~\cite{CEG}. Similarly, the smallest known uncolorable vector set in dimension~3, due to Conway and Kochen~\cite[pp.~114, 197]{Peres:quantum}, consists of integer vectors on the cube with vertices $(\pm 2, \pm 2, \pm 2)$.
Another uncolorable set of~37 integer vectors in dimension~3 was produced by Bub~\cite{Bub:Schutte}. More recently, large numbers of integer KS vector sets were produced computationally by Pavi\v{c}i\'{c}~\cite{Pavicic:exhaustive} and Pavi\v{c}i\'{c}-Megill~\cite{PavicicMegill:automated}.

In this section we construct KS-uncolorable sets of integer vectors that have noteworthy number-theoretic features. The important feature of these sets of vectors are described in terms of a numerical value associated to the sets, which we now define. 
Let $q \colon \Z^d \to \Z_{\geq 0}$ denote the Euclidean quadratic form, defined on a column vector $v = (v_1,\dots,v_d)^T$ by
\begin{align*}
q(v) &= \|v\|^2 = v^T v 
= v_1^2 + \cdots + v_d^2.
\end{align*}
If $S \subseteq \Z^d$ is a Kochen-Specker uncolorable set of integer vectors, we will be interested in the numerical invariant
\[
N(S) := \operatorname{lcm}\{q(v) \mid v \in S\}.
\]
If $v$ is an integer column vector, then the orthogonal projection onto the line through $v$ is given by the matrix 
\[
P_v = q(v)^{-1} \, v \, v^T,
\]
whose entries are evidently \emph{rational} and lie in the subring $\Z[1/q(v)] \subseteq \Q$. 
Thus for all $v \in S$, the projection matrix $P_v$ has entries in $\Z[1/N(S)]$. 

For instance, if $S_1$ denotes any of the sets of 4-dimensional integer vectors in~\cite{CEG, Kernaghan, Peres:KS}, 
then satisfies $q(S_1) = \{1,2,4\}$ and thus $N(S_1) = 4$. 
Similarly, for the set $S_2 \subseteq \Z^3$ constructed in~\cite{Bub:Schutte}, we have $q(S_2) =  \{1,2,3,5,6,30\}$ so that $N(S_2) = 30$. The Conway-Kochen set $S_3 \subseteq \Z^3$ is traditionally represented~\cite[p.~114]{Peres:quantum}  by vectors on the cube with vertices $(\pm 2, \pm 2, \pm 2)$, so that it satisfies $q(S_3) = \{4, 5, 6, 8, 9, 12\}$ and $N(S_3) = 360$.
However, some of these vectors such as $(2,2,2) = 2 \cdot (1,1,1)$ or $(2,0,0) = 2 \cdot (1,0,0)$ can be replaced by a smaller collinear integer vector to obtain an equivalently uncolorable set $S_4 \subseteq \Z^3$ that satisfies $q(\setS_4) = \{1, 2, 3, 5, 6, 9\}$ and $N(S_4) = 90$.

\separate 

Given a fixed positive integer $N$, we may take the reverse perspective and ask: \emph{is there an uncolorable set of integer vectors whose corresponding projection matrices have entries in $\Z[1/N]$?}
Given the prime factorization $N = \prod p_i^{e_i}$, recall that the \emph{radical} is defined as $\rad N = \prod p_i$. Noting that the rings $\Z[1/N] = \Z[1/\rad N]$ coincide, to answer this question it suffices to consider only squarefree values of $N$.

For positive integers $N$ and $d$, we define
\begin{align*}
\setS_d(N) &= \{ v \in \Z^d : q(v) \mbox{ is a unit in } \Z[1/N]\} \\
 &= \{ v \in \Z^d : q(v) \mbox{ divides a power of } N\}.
\end{align*}
Note that $\setS_d(N)$ is also the set of those $v \in \Z^d$ such that the primes occurring in the factorization of $q(v)$ form a subset of those occurring in $N$, and that $\setS_d(N) = \setS_d(\rad N)$. As above, if $v \in \setS_d(N)$ then its corresponding projection matrix satisfies $P_v \in \M_d(\Z[1/N])$.

We wish to address the following: 

\begin{question}
\label{q:colorable}
For which positive squarefree integers $N$ and dimensions $d \geq 3$ is the set of vectors $\setS_d(N)$ Kochen-Specker uncolorable?
\end{question}

In contrast to much of the existing literature on KS uncolorable vector sets, we are \emph{not} interested in isolating minimal uncolorable sets. Rather, we are interested in the number-theoretic question of which integers $N$ give rise to uncolorable or colorable sets in a given dimension. For instance, given the sets $S_2$ of Bub and $S_3$ or $S_4$ of Conway-Kochen discussed above, we have $S_2, S_3, S_4 \subseteq \setS_3(30)$, so each of these vector sets independently shows that $\setS_3(30)$ has no KS coloring. Although $S_3$ and $S_4$ have fewer vectors, from the perspective of this number-theoretic measure we will prefer to work with the set $S_2$ because $N(S_2) = 30$ while $N(S_4) = 90$.

It is clear that if an integer $M$ divides $N$, then $\setS_d(M) \subseteq \setS_d(N)$. It follows that for $M \mid N$,
\begin{align*}
\setS_d(N) \mbox{ colorable} &\implies \setS_d(M) \mbox{ colorable,} \\
\setS_d(M) \mbox{ uncolorable} &\implies \setS_d(N) \mbox{ uncolorable.}
\end{align*}
We will use these implications freely below.

\subsection{Uncolorability in dimension~3}

We first focus on the classical case of dimension $d = 3$. 
It was shown in~\cite[Theorem~3.4]{BMR} that if $N$ is not divisible by~2 or not divisible by~3, then $\setS_3(N)$ has a KS coloring. 
As an immediate consequence,
\begin{center}
if $\setS_3(N)$ is KS uncolorable, then $6$ divides $N$.
\end{center}
On the other hand, Bub's uncolorable set $S \subseteq \Z^3$ with $N(S) = 30$ implies that
\begin{center}
if $30$ divides $N$, then set $\setS_3(N)$ is KS uncolorable.
\end{center}

In light of the above facts, one might naturally guess that $\setS_3(6)$ has no Kochen-Specker coloring. However, to date we have not managed to find a proof or refutation of this statement.
Lacking such an answer, we wondered whether there exists an integer~$N$, divisible by~$6$ but not~$5$, such that $\setS_3(N)$ is KS uncolorable.  In this subsection we exhibit just such an example, for the specific value
\[
N = 2 \cdot 3 \cdot 7 \cdot 11 = 462.
\]

This will be proved by constructing an uncolorable subset $Q \subseteq \setS_3(462)$, which we now describe.
When considering Kochen-Specker colorings of sets of vectors, each vector is intended to represent the rank-one projection onto the line through that vector. For this reason it suffices to consider sets of pairwise noncollinear vectors. When restricting attention to integer vectors, this means that we may consider only those vectors whose entries have greatest common divisor equal to~1. We will call such integer vectors \emph{primitive}.

Furthermore, for each nonzero vector $v$, our vector sets need to only contain one of the two vectors $\{v, -v\}$. For convenience, we will say that a vector $v = (v_1,v_2,v_3) \in \Z^3 \setminus \{0\}$ is \emph{well-signed} if either: 
\begin{itemize}
\item $v$ has only one nonzero entry which is positive,
\item $v$ has two nonzero entries and its first nonzero entry is positive, or
\item $v$ has three nonzero entries, at least two of which are positive.
\end{itemize}
For instance, the vectors $(1,0,0)$, $(0,1,-1)$, $(1,1,1)$, and $(1,-1,1)$ are well-signed while $(-1,0,0)$, $(0,-1,1)$, $(-1,-1,-1)$, and $(-1,1,-1)$ are not. It is clear that for each $v \in \Z^3 \setminus \{0\}$, exactly one of $v$ or $-v$ is well-signed. It follows that if a set $S \subseteq \Z^3$ consists of primitive, well-signed vectors, then the vectors in $S$ must be noncollinear.

For $n = 1,2,3,6,21,33,77$, we will define subsets $Q_n \subseteq \setS_3(462)$ such that every $v \in Q_n$ satisfies $q(v) = n$. 
For the values $n = 1,2,3,6,21$, we let $Q_n$ denote the set of all well-signed primitive integer vectors $v$ such that $q(v) = n$. It follows that:
\begin{itemize}
\item $Q_1 = \{(1,0,0),\ (0,1,0),\ (0,0,1)\}$.
\item $Q_2 = \{(1,1,0),\ (1,0,1),\ (0,1,1),\ (1,-1,0),\ (1,0,-1), (0,1,-1)\}$.
\item $Q_3 = \{(1,1,1),\ (1,1,-1),\ (1,-1,1),\ (-1,1,1)\}$.
\item $Q_6$ contains $12$ such vectors, and they are the well-signed vectors whose entries up to sign are $1,1,2$.
\item $Q_{21}$ contains $24$ vectors, and they are the well-signed vectors whose entries up to sign are $1,2,4$.
\end{itemize}
For the remaining values $n=33,77$, we define the sets as follows.
\begin{itemize}
\item $Q_{33}$ is the set of all well-signed vectors whose entries up to sign are $2,2,5$. There are~$12$ vectors in this set. 
\item $Q_{77}$ is the set of all well-signed vectors whose entries up to sign are $2,3,8$. There are~$24$ vectors in this set. 
\end{itemize}
As in the previouis cases, the vectors $v \in Q_n$ are all well-signed and primitive, satisfying $q(v) = n$. However, there are certain vectors with the ``correct'' value of $q(v)$ that are excluded from these sets, such as $(1,4,4) \notin Q_{33}$ and $(4,5,6) \notin Q_{77}$.

Finally, we will define $Q \subseteq \setS_3(462)$ to be the disjoint union
\begin{equation}\tag{$\star$} \label{eq:Q}
Q = Q_1 \cup Q_2 \cup Q_3 \cup Q_6 \cup Q_{21} \cup Q_{33} \cup Q_{77}.
\end{equation}
This is a set of~$85$ integer vectors that are primitive, well-signed, and thus noncollinear. 
This set and its uncolorability were discovered using Mathematica software.\footnote{The code is publicly available at \url{https://github.com/manny-reyes/Kochen_Specker_Colorability/tree/main/462_vector_set}.}
However, we are able to produce a human-readable proof of this fact as follows.

\begin{theorem}
\label{thm:Q}
There is no Kochen-Specker coloring of $Q \subseteq \setS_3(462)$.
\end{theorem}

\begin{proof}
Assume toward a contradiction that $Q$ has a KS coloring. Because $Q$ is invariant under permutation of coordinates, we may assume without loss of generality that in the orthogonal triple
\[
\{(1,0,0),\ (0,1,0),\ (0,0,1)\}
\]
the vector $(1,0,0)$ is colored~$1$ while the other two are colored~$0$. It follows that in the triple
\[
\{(1,0,0),\ (0,1,1),\ (0,1,-1)\}
\]
the second and third vectors must be colored~0. We also see that in the orthogonal triple
\[
\{(0,1,0),\ (1,0,1),\ (1,0,-1)\}
\]
either the second or third vector must be colored~1. Because multiplication by the diagonal reflection matrix $R_z = \diag(1,1,-1)$ leaves the coordinate axes invariant while interchanging the vectors $(1,0,\pm1)$, we may again assume without loss of generality that $(1,0,1)$ is colored~0 and $(1,0,-1)$ is colored~1. Finally, in the orthogonal triple
\[
\{(0,0,1),\ (1,1,0),\ (1,-1,0)\}
\]
either the second or third vector must be colored~1. Since multiplication by the reflection matrix $R_y = \diag(1,-1,1)$ again leaves the coordinate axes and the lines through $(1,0,\pm1)$ invariant while interchanging the vectors $(1,\pm1,0)$, we may assume without loss of generality that $(1,1,0)$ is assigned~0 and $(1,-1,0)$ is~1. 

To summarize, we are assuming without loss of generality the following assignments of colors:
\begin{itemize}
\item \textbf{Color~1:} $(1,0,0)$, $(1,0,-1)$, $(1,-1,0)$
\item \textbf{Color~0:} $(0,1,0)$, $(0,0,1)$, $(0,1,1)$, $(0,1,-1)$, $(1,0,1)$, $(1,1,0)$.
\end{itemize}
We now deduce the following sequence of colorings from the listed orthogonal triples:
\begin{align*}
\{(1,0,-1),\ (1,-1,1)\} \mbox{ orthogonal} &\implies (1,-1,1) \mapsto 0, \\ 
\{(0,1,1),\ (1,-1,1),\ (2,1,-1)\} \mbox{ orthogonal} &\implies (2,1,-1) \mapsto 1, \\ 
\{(2,1,-1),(-3,8,2)\} \mbox{ orthogonal} &\implies (-3,8,2) \mapsto 0. 
\end{align*}
In a similar manner, we have the following colorings:
\begin{align*}
\{(1,-1,0),\ (1,1,-1)\} \mbox{ orthogonal} &\implies (1,1,-1) \mapsto 0, \\ 
\{(1,0,1),\ (1,1,-1),\ (-1,2,1)\} \mbox{ orthogonal} &\implies (-1,2,1) \mapsto 1, \\ 
\{(-1,2,1),\ (4,1,2)\} \mbox{ orthogonal} &\implies (4,1,2) \mapsto 0, \\ 
\{(1,-1,0),\ (2,2,-5)\} \mbox{ orthogonal} &\implies (2,2,-5) \mapsto 0, \\ 
\{(4,1,2),\ (2,2,-5),\ (-3,8,2)\} \mbox{ orthogonal} &\implies (-3,8,2) \mapsto 1. 
\end{align*}
But this contradicts the previous deduction that $(-3,8,2)$ must be colored~0, completing the proof.
\end{proof}

\subsection{A computational method for detecting uncolorability}
\label{sec:computer}

Our further explorations into the uncolorability of $\setS_d(N)$ were aided by computer searches. This consists of generating the vectors $v \in \setS_d(N)$ such that $q(v)$ divides $N^e$ for some power $e \geq 1$, and testing whether the set of vectors is colorable or not.
In order to test colorability of the vectors, we have found it useful to apply the method of~\cite[Construction~3.1.3]{Salt} to convert the question of whether a Kochen-Specker coloring exists into an integer linear programming problem. We describe that method below.

Fix a finite set of nonzero vectors $S = \{v_j\}_{j = 1}^s \subseteq \C^d$.
Let $G_S$ denote the orthogonality graph of $S$, whose vertices are the vectors in $S$ and where there is an edge connecting $v,w \in S$ if and only if $v \perp w$.
Let $\{C_i\}_{i=1}^r$ be an enumeration of the maximal cliques of the graph $G_S$; in other words, these are the maximal subsets of $V$ consisting of pairwise orthogonal vectors. We necessarily have $|C_i| \leq d$ since orthogonal sets of nonzero vectors are linearly independent.

Define an $r \times s$ matrix $M = (m_{ij})$ by setting
\[
m_{ij} = \begin{cases}
1, & v_j \in C_i, \\
0, & v_j \notin C_i.
\end{cases}
\]
Let $\overrightharp{1}$ be the $r$-dimensional vector with all entries~1, and define $b = (b_i)_{i=1}^r$ by
\[
b_i = \begin{cases}
1 & |C_i| = d, \\
0 & |C_i| < d.
\end{cases}
\]
Then a Kochen-Specker coloring of $V$ is equivalent to a solution $(x_j) \in \{0,1\}^s$ to the integer linear programming problem (with zero constraint function)
\begin{align*}
b \leq Mx \leq \overrightharp{1}.
\end{align*}
Indeed, values $x_j$ correspond to the color of the vector $v_j$.
The $i$th entry of the matrix product $\sum_{j=1}^s m_{ij} x_j$ counts the number of vectors in $C_i$ that are colored~1. 
The upper bound ensures that each clique has at most one vector assigned value~1, while the lower bound ensures that each orthogonal basis in $V$ has exactly one vector assigned value~1. 

\separate

We implemented this method using Python code.\footnote{The code can be found in the folder \url{https://github.com/manny-reyes/Kochen_Specker_Colorability/tree/main/Dimension_d_vector_coloring}.} For a positive integer $N$ and dimension $3 \leq d \leq 6$, the algorithm proceeds as follows:
\begin{enumerate}
\item Produce the set $S$ all primitive integer vector solutions $x = (x_1, \dots, x_d)$ to $\sum_{i=1}^d x_i^2 = m$ for every $m \mid N$. Remove from $S$ any solution $x$ whose first nonzero entry is negative, to guarantee that no two vectors are collinear. (This uses a diophantine equation solver from the package \texttt{SymPy}.)
\item Produce the orthogonality graph $G_S$ of this set of vectors by computing all dot products of pairs in $S$. (The graph is a class in the \texttt{NetworkX} package.)
\item Enumerate all maximal cliques of $G_S$. (This is a built-in function of \texttt{NetworkX}.)
\item Produce the matrix $M$ and lower bound vector $b$ defined above.
\item Check whether the integer linear programming problem $b \leq Mx \leq \overrightharp{1}$ has a solution. (This uses the mixed integer linear programming function of \texttt{SciPy}.)
\end{enumerate}

Using this code we have been able to explore whether $\setS_d(N)$ is colorable for various values of $N$ across dimensions $3 \leq d \leq 6$. 
Some of our computational findings in dimension $d = 3$ are summarized in Table~\ref{tab:dim 3}. Our computations verify those of~\cite[\S 3.3.1]{Salt} for all overlapping values of $N$, including the interesting fact that $\setS_3(714)$ is uncolorable.

\begin{table}
\caption{Colorability of vectors $v \in \setS_3(N)$ with $q(v) \leq N^e$, where $p$ denotes a prime.}
\begin{tabular}{|c |c| c|}
\hline
$N$ & $e$ & result \\
\hline \hline
6 & 10 & colorable \\
\hline 
$30 = 6 \cdot 5$ & 1 & uncolorable \\
\hline
\begin{tabular}{c}
$6 \cdot p$, \\
$p \in \{7, 11, 13,\dots , 103\}$
\end{tabular} & 2 & colorable \\
\hline
$462 = 6 \cdot 7 \cdot 11$ & 1 & uncolorable \\
\hline
$714 = 6 \cdot 7 \cdot 17$ & 1 & uncolorable \\
\hline 
\begin{tabular}{c}
$6 \cdot 7 \cdot p$, \\
$p \in \{13, 19, 23, \dots, 103\}$
\end{tabular} & 1 & colorable \\
\hline
\end{tabular}
\label{tab:dim 3}
\end{table}

\separate

This method of testing colorability can be adapted in a straightforward way from colorings of vectors to colorings of projection matrices over, say, finite fields (as defined in Section~\ref{sec:rings} below). Let $S$ be the set of projection matrices in question, with the zero matrix removed for convenience. One then produces the orthogonality graph $G_S$ for that set of matrices. One then finds the maximal cliques of $G_S$ and defines a similar matrix $M$. With the zero matrix removed, the maximal cliques are pairwise orthogonal sets of projections that sum to the identity matrix. 
In this case the lower and upper bound vectors both equal $\overrightharp{1}$, and we are in fact seeking a solution $x \in \{0,1\}^r$ to $Mx = \overrightharp{1}$.

\subsection{Uncolorability in dimensions $d \geq 4$}

In this subsection we will address Question~\ref{q:colorable} in higher dimensions. A natural first question is whether there is any relationship between KS colorability in a given dimension and that in higher dimensions. The following lemma tells us that uncolorability passes directly to higher dimensions.

\begin{lemma}\label{lem:higher dim}
Let $N$ and $d$ be positive integers. If $\setS_{d}(N)$ is KS-uncolorable then $\setS_{d+k}(N)$ is also KS-uncolorable for all integers $k \geq 0$. 
\end{lemma}

\begin{proof}
The method of proof is similar to that of~\cite[Lemma~2.17]{BMR}. Assume that $\setS_d(N)$ is KS-uncolorable; proceeding by induction, it suffices to show that $\setS_{d+1}(N)$ is also KS-uncolorable. 

Assume toward a contradiction that there is a KS coloring $c \colon \setS_{d+1}(N) \to \{0,1\}$. Then one of the standard basis vectors $e_1, \dots, e_{d+1}$ is colored~1 while the rest are colored~0. Assume without loss of generality that $c(e_{d+1}) = 0$.
Let $U \colon \C^d \hookrightarrow \C^{d+1}$ be the isometry $U(v) = (v,0)$ that appends zero as the last entry. This preserves inner products, and it restricts to a map $U \colon \Z^d \hookrightarrow \Z^{d+1}$ that preserves dot products. 
In particular, $q(U(v)) = q(v)$ for $v \in \Z^d$. Thus we have an injective map
\[
U \colon \setS_d(N) \hookrightarrow \setS_{d+1}(N),
\]
that perserves orthogonality. 

This induces a coloring $c \circ U \colon \setS_d(N) \to \{0,1\}$, and properties (KS1) and (KS2) are obviously inherited from the coloring $c$. To see that axiom (KS3) is also satisfied, suppose that $v_1, \dots, v_d \in \setS_d(N)$ are orthogonal. Then $U(v_1), \dots, U(v_d), e_{d+1} \in \setS_{d+1}(N)$ are also orthogonal. Because $c(e_{d+1}) = 0$, we must have $c(U(v_i)) = 1$ for some $i$. Thus $c \circ U$ is a Kochen-Specker coloring of $\setS_d(N)$, contradicting the hypothesis.
\end{proof}

We have used the method of Subsection~\ref{sec:computer} to investigate colorability of $\setS_d(N)$ for the range of dimensions $d = 4, 5, 6$. Some of these findings are summarized in Table~\ref{tab:higher dim}. The most notable discovery is the uncolorability of $\setS_6(3)$, which we record in the theorem below. The uncolorable set consists of those vectors $v \in \Z^6$ such that $q(v) \in \{1,3\}$. Choosing a single representative out of each pair $\pm v$ of such vectors on the same line yields a set of 86 vectors. (This includes the six standard basis vectors, and $\binom{6}{3} \cdot 2^2 = 80$ other vectors with three nonzero entries that are $\pm 1$ with, say, their first entries positive.) We have not searched to find the smallest possible uncolorable subset of these vectors.

\begin{table}
\caption{Colorability of vectors $v \in \setS_d(N)$ with $q(v) \leq N^e$.}
\begin{tabular}{|c| c | c | c|}
\hline
$d$ & $N$ & $e$ & result \\
\hline \hline
4 & 2 & 2 & uncolorable \\
\hline 
4 & $p \in \{3, 5, \dots, 13 \}$ & 3 & colorable \\
\hline
5 & $p \in \{5, 7, 11, 13\}$ & 2 & colorable \\
\hline
6 & 3 & 1 & uncolorable \\
\hline 
6 & $p \in \{5, 7, 11,\dots , 19\}$ & 1 & colorable \\
\hline
\end{tabular}
\label{tab:higher dim}
\end{table}

\begin{theorem}\label{thm:dim 4 and 6}
Let $d$ be a positive integer.
\begin{enumerate}
\item For $d \geq 4$, the set $\setS_d(2)$ has no Kochen-Specker coloring.
\item For $d \geq 6$, the set $\setS_d(3)$ has no Kochen-Specker coloring.
\end{enumerate}
\end{theorem}

\begin{proof}
(1) As noted earlier, there is a KS uncolorable set $S \subseteq \Z^4$ constructed in~\cite{Peres:KS} such that $N(S) = 4$. Thus $\setS_4(2) \supseteq S$ is also uncolorable. The case $d \geq 4$ now follows from Lemma~\ref{lem:higher dim}.

(2) As shown in Table~\ref{tab:higher dim}, the set of vectors $v \in \setS_6(3)$ satisfying $q(v) \in \{1,3\}$  is KS uncolorable. So the result follows from Lemma~\ref{lem:higher dim} again.
\end{proof}

The dimension $d = 4,5$ results in Table~\ref{tab:higher dim} suggested to us the possibility that $\setS_d(N)$ might be KS-colorable for odd $N$ in these dimensions. This led us to apply the method described Subsection~\ref{sec:computer} to investigate whether the projection matrices over the field $\F_2$ with two elements were Kochen-Specker colorable.\footnote{This code is available at \url{https://github.com/manny-reyes/Kochen_Specker_Colorability/tree/main/Projection_coloring}.} Our findings are recorded in Table~\ref{tab:projections}.

\begin{table}
\caption{Colorability of $\Proj(\M_d(\F_2))$.}
\begin{tabular}{|c|c|}
\hline
$d$ & result \\
\hline \hline
3 & colorable \\
\hline 
4 & colorable \\
\hline
5 & uncolorable \\
\hline
\end{tabular}
\label{tab:projections}
\end{table}

Colorability in dimension~4 verifies this guess, but the approach via projection matrices was complicated by the uncolorability in dimension~5. 
Nevertheless, we were able to verify our initial guess with an explicit proof in both of dimensions~4 and~5 as follows.

\begin{proposition}\label{prop:dim 4 and 5}
If $N$ is an odd integer and $d = 4,5$, then $\setS_d(N)$ has a Kochen-Specker coloring. In fact, there is a coloring of the set $O_d = \bigcup_{2 \nmid N} \setS_d(N)$. 
\end{proposition}

Before proving the result, we set up some notation. Below we consider the ordinary dot product of vectors in the $\F_2$-vector spaces $\F_2^d$. We define the following set of $d$-dimensional vectors over the field with two elements:
\begin{align*}
\setS_d(\F_2) &= \{v \in \F_2^d \mid v \cdot v \neq 0\} \\
&= \{v \in \F_2^d \mid v \mbox{ has an odd number of nonzero entries}\}.
\end{align*}
Denote the canonical surjection $\pi \colon \Z^d \twoheadrightarrow (\Z/2\Z)^d = \F_2^d$ by $v \mapsto \overline{v}$; with a slight abuse of notation we use the same notation for all $d \geq 1$. Note that this is compatible with the dot product of vectors in the sense that for all $v, w \in \Z^d$,
\[
\overline{v} \cdot \overline{w} = \overline{v \cdot w}.
\]
In particular, this means that $\pi$ preserves orthogonality of vectors. Furthermore, if $q(v) = v \cdot v$ is odd, then $\pi(v) \in \setS_d(\F_2)$. This means that $\pi$ restricts to an orthogonality-preserving map 
\[
\pi \colon O_d = \bigcup_{2 \nmid N} \setS_d(N) \to \setS_d(\F_2).
\]

\begin{proof}[Proof of Proposition~\ref{prop:dim 4 and 5}]
We let $e_1, \dots, e_d \in \F_2^d$ denote the standard basis vectors below.  

\textbf{Case $d = 4$:} For $i = 1,\dots,4$ define vectors $f_i = (1,1,1,1) - e_i$, so that $f_i$ has a zero in the $i$th component and ones in all other components. Then we have
\[
\setS_4(\F_2) = \{e_1, e_2, e_3, e_4, f_1, f_2, f_3, f_4\},
\]
One can easily verify that $e_i$ is orthogonal to $f_j$ if and only if $i = j$, and that the only orthogonal quadruples of vectors in this set are $\{e_1, e_2, e_3, e_4\}$ and $\{f_1, f_2, f_3, f_4\}$. The coloring $c \colon \setS_4(\F_2) \to \{0,1\}$ that sends $e_1, f_2 \mapsto 1$ and all other vectors to zero clearly satisfies the conditions (KS2)--(KS3) of a Kochen-Specker coloring.

Now if we compose the orthogonality-preserving map above with this coloring, we obtain a coloring
\[
c \circ \pi \colon O_4 \to \setS_4(\F_2) \to \{0,1\}.
\]
Because $\pi$ preserves orthogonality, axioms (KS2) and (KS3) still hold for this coloring. The fact that (KS1) still holds follows from the fact that collinear vectors in $\setS_d(N)$ must be of the form $v$ and $\lambda v$ where $\lambda$ is an integer dividing a power of $N$; this means $\lambda$ is odd, so that $\pi(\lambda v) = \pi(v)$ and both vectors are assigned the same color. 
Thus this is a Kochen-Specker coloring of $O_4$.

\textbf{Case $d = 5$:} 
Denote $g = (1,1,1,1,1)$, and for any pair of integers $1 \leq i < j \leq 5$ we define the vector $f_{ij} = g -e_i - e_j$ to have zeros in the $i$th and $j$th coordinates and ones in the other three coordinates. For notational covenience, we also denote $f_{ji} = f_{ij}$. This time we have
\[
\setS_5(\F_2) = \{e_i \mid i = 1, \dots 5\} \cup \{f_{ij} \mid 1 \leq i < j \leq 5\} \cup \{g\}.
\]
There are 6 orthogonal quintuples in this family. One is the standard basis $\{e_1, \dots, e_5\}$, and the others are of the form $\{e_i\} \cup \{f_{ij} \mid j \neq i\}$ for each $i = 1, \dots, 5$. For instance, one such orthogonal set is $\{e_1, f_{12}, f_{13}, f_{14}, f_{15}\}$. The coloring $c \colon \setS_5(\F_2) \to \{0,1\}$ that sends $e_1, f_{23}, f_{45} \mapsto 1$ and all other vectors to~0 satisfies (KS2) and (KS3), so that as before we may compose $c \circ \pi \colon O_5 \to \{0,1\}$ to obtain a Kochen-Specker coloring.
\end{proof}

This has the curious consequence that Question~\ref{q:colorable} can be fully answered in dimensions~4 and~5, a situation which seems unique to those dimensions.

\begin{corollary}\label{cor:dim 4 and 5}
For $d = 4,5$ and any positive integer $N$, the set $\setS_d(N)$ is Kochen-Specker uncolorable if and only if $N$ is even.
\end{corollary}

\section{Consequences for symmetric matrices over commutative rings}
\label{sec:rings}

We now apply the results of Section~\ref{sec:uncolorable} to the study of contextuality in the purely algebraic setting, culminating in a proof of Theorem~\ref{thm:map to commutative}. 
The link from vector colorings to algebraic hidden states occurs in a multi-step process: uncolorability of vectors implies uncolorability of projection matrices, which in turn implies that a certain spectrum is empty (Lemma~\ref{lem:coloring}), which finally implies the non-existence of morphisms to commutative rings. 
To explain these connnections, we briefly provide several definitions from~\cite[Section~2]{KochenSpecker} and~\cite[Section~2]{BMR}, to which readers are referred for more details.

A \emph{partial Boolean algebra} is a set $B$ equipped with a reflexive, symmetric binary operation $\odot$ of \emph{commesaurability}, distinguished elements $0,1 \in B$, a unary operation $\neg$ of negation, and partially defined binary operations of meet $\wedge$ and join $\vee$ for commeasurable pairs of elements, such that every pairwise commeasurable set $S \subseteq B$ is contained in a pairwise commeasurable set $C \subseteq B$ containing $0$ and $1$ for which the restricted negation, meet, and join make $C$ into an ordinary (``total'') Boolean algebra. 
The colorings of Kochen and Specker correspond to \emph{homomorphisms} of partial Boolean algebras $B \to \mathbf{2}$, where $\mathbf{2} = \{0,1\}$ is the ordinary two-element Boolean algebra. We simply refer to such morphisms as Kochen-Specker colorings of a partial Boolean algebra as in~\cite[Theorem~2.13]{BMR}.

If $K$ is a commutative ring, then the symmetric (under the transpose operation) idempotent elements of the matrix ring $\M_n(K)$ will be called \emph{projections} by analogy with the case of real matrices. The set of all projections $\Proj(\M_n(K))$ forms a partial Boolean algebra, where commeasurability is given by commutativity in $\M_n(K)$, negation is the complement $\neg \, e = I - e$, and meet and join are given by $e \wedge f = ef$ and $e \vee f = e + f - ef$. 
For $K = \R$, this yields the usual partial Boolean algebra of projection matrices. (A suitable modification to the case where $K$ is equipped with an involution allows one to recover the partial Boolean algebra of complex projection matrices, as well.)

\separate

Bridging the gap between colorability of projections and morphisms of partial rings requires prime partial ideals, as defined in~\cite{Reyes}. Let $R$ be a partial ring.
A subset $I \subseteq R$ is a \emph{partial ideal} if it satisfies the following conditions for all $a, b \in R$:
\begin{itemize}
\item $a,b \in I$ and $a \odot b \implies a + b \in I$,
\item $b \in I$ and $a \odot b \implies ab \in I$.
\end{itemize}
A partial ideal $P \subseteq R$ is \emph{prime} if it additionally satisfies $1 \notin P$ and, for all $a,b \in R$, 
\begin{itemize}
\item $a \odot b$ and $ab \in P \implies a \in P$ or $b \in P$.
\end{itemize}
It is straightforward to see that a subset $P \subseteq R$ is a (prime) partial ideal if and only if, for every commeasurable total subring $C \subseteq R$, the intersection $C \cap P$ is a (prime) ideal of $C$. 

The \emph{partial spectrum} $\pSpec(R)$ is defined to be the set of all prime partial ideals of $R$. 
This forms a contravariant functor from the category of partial rings to the category of sets in the usual way, where a morphism $f \colon R \to S$ induces a map $\pSpec(S) \to \pSpec(R)$ by $P \mapsto f^{-1}(P)$.
The link between the partial spectrum of $R$ and Kochen-Specker colorings is given by the following observation from~\cite[(2.15)]{BMR} . If there exists a prime partial ideal $P \in \pSpec(R)$, then we obtain a Kochen-Specker coloring on the partial Boolean algebra $B = \Idpt(R)$ of idempotents in $R$, given by the following assignment for $e \in B$: 
\[
e \mapsto 
\begin{cases}
0, & e \in B \cap P, \\
1, & e \in B \setminus P.
\end{cases}
\]
On the other hand, if $\pSpec(R) = \varnothing$ then there is no morphism of partial rings $R \to C$ for any nonzero commutative (total) ring $C$ by~\cite[Lemma~2.4]{BMR}.
In this way, the non-existence of prime partial ideals can be viewed as an algebraic manifestation of contextuality.

Note that if $K$ is a commutative ring and $R = \M_n(K)_\sym$ is a partial algebra of symmetric matrices, then the partial Boolean algebra of idempotents recovers the partial Boolean algebra of projections discussed above: $\Idpt(R) = \Proj(\M_n(K))$. Thus uncolorability of $\Proj(\M_n(K))$ implies that there is no morphism of partial rings $\M_n(K)_\sym \to C$ for any nonzero commutative ring $C$.

We record the following sequence of implications relating KS uncolorability to an empty partial spectrum.

\begin{lemma}\label{lem:coloring}
For positive integers $N$ and $d$, consider the following statements:
\begin{enumerate}[label=\textnormal{(\roman*)}]
\item $\setS_d(N)$ has no Kochen-Specker coloring; 
\item $\Proj(\M_d(\Z[1/N]))$ has no Kochen-Specker coloring;
\item $\M_d(\Z[1/N])_\sym$ has no prime partial ideals.
\end{enumerate}
Then $\mathrm{(i)} \implies \mathrm{(ii)} \implies \mathrm{(iii)}$.
\end{lemma}

\begin{proof}
(i)$\implies$(ii): Each vector in $\setS_d(N)$ defines a rank-1 projection in $\M_d(\Z[1/N])$, with orthogonality of vectors corresponds to orthogonality of the associated projections. Thus every Kochen-Specker coloring $\Proj(\M_d(\Z[1/N])) \to \{0,1\}$ induces a corresponding KS coloring of $\setS_d(N)$. If $\setS_d(N)$ is KS uncolorable, this means that the projections of $\M_d(\Z[1/N])$ are also KS uncolorable.

(ii)$\implies$(iii) is proved in~\cite[Corollary~2.16]{BMR}.
\end{proof}

\separate

The next lemma explains how uncolorability of $\setS_d(M_1)$ and $\setS_d(M_2)$ can be used to deduce information about $\M_d(\Z[1/N])_\sym$ where $N = \gcd(M_1, M_2)$. The statement is phrased a bit more generally.

\begin{lemma}\label{lem:Bezout}
For an integer $d \geq 3$, suppose that $M_1, \dots, M_r$ are positive integers such that each $\pSpec(\M_d(\Z[1/M_i])_\sym) = \varnothing$. Then for any integer $N$ divisible by $M = \gcd(M_1, \dots, M_r)$, the partial ring $\M_d(\Z[1/N])_\sym$ has no prime partial ideals.
\end{lemma}

\begin{proof}
Assume toward a contradiction that there exists a prime partial ideal $P' \in \pSpec(\M_d(\Z[1/N])_\sym)$. 
Restricting to the partial subring $\M_d(\Z[1/M])_\sym \subseteq \M_d(\Z[1/N])_\sym$, we obtain a prime partial ideal $P = P' \cap \M_d(\Z[1/M])_\sym$. 
Fix $1 \leq i \leq r$, and note that $\Z[1/M] \subseteq \Z[1/M_i]$ because $M$ divides $M_i$. We claim that 
\[
M_i I_d \in P.
\]
Denote $Q = \left\{\frac{1}{M_i^e} \cdot x \mid x \in P, \ e \geq 0\right\} \subseteq \M_d(\Z[1/M_i])_\sym$. It is straightforward to verify that $Q$ is also a partial ideal of the partial ring $R = \M_d(\Z[1/M_i])_\sym$, and that it still satisfies the condition that if $r,s \in R$ are commeasurable and $rs \in Q$, then $r$ or $s$ lies in $Q$. But $\pSpec(R) = \varnothing$, so $Q$ cannot be a prime partial ideal. This is only possible if $Q = R$. Thus $I_d \in Q$, from which it follows that $(M_i I_d)^e = M_i^e \cdot I_d \in P$ for some $i \geq 0$. Because $P$ is prime, it follows that $M_i I_d \in P$ as claimed.

Finally, use B\'{e}zout's identity to write $M = \sum c_i M_i$ for some $c_i \in \Z$. Because the elements $M_i I \in \M_d(\Z[1/M])_\sym$ are pairwise commeasurable and $P$ is a partial ideal, we must have
\[
M I_d = c_1(M_1I_d) + \cdots + c_r (M_r I_d) \in P.
\]
But then $I_d = (\frac{1}{M} I_d)(MI_d) \in P$, contradicting that $P$ (and thereby $P'$) is prime. 
\end{proof}

\separate

Using the strategies described above, can apply the results of Section~\ref{sec:uncolorable} to characterize exactly which of the rings of the form $K = Z[1/N]$ (i.e., finitely generated subrings of $\Q$) have $\pSpec(\M_3(K)_\sym) = \varnothing$. We also have partial information in dimensions~4 and~5, and the strongest possible result holding for dimensions $d \geq 6$.

\begin{theorem}\label{thm:prime}
Let $N$ and $d$ denote a positive integers below.
\begin{enumerate}
\item The partial ring $\M_3(\Z[1/N])_\sym$ has a prime partial ideal if and only if~$N$ is not divisible by~$6$.
\item For $d = 4,5$, if $N$ is even then $\pSpec(\M_d(\Z[1/N])) = \varnothing$.
\item If $d \geq 6$, then $\pSpec(\M_d(\Z)) = \varnothing$.
\end{enumerate}
\end{theorem}

\begin{proof}
(1) 
If $N$ is not divisible by~$6$, then it is relatively prime to~2 or~3 (possibly both). This means that $N$ is a unit modulo~2 or~3, so that we have a ring homomorphism from $\Z[1/N] \to \F_p$, and consequently a morphism $\M_3(\Z[1/N])_\sym \to \M_3(\F_p)_\sym$, for either $p = 2,3$. By~\cite[Theorem~3.4]{BMR} there exists a prime partial ideal of $\M_3(\F_p)_\sym$ for each of these values of~$p$, and its preimage gives a prime partial ideal of $\M_3(\Z[1/N])_\sym$. 

Now suppose that $N$ is divisible by~6. By~\cite{Bub:Schutte}, the set $\setS_3(30)$ has no Kochen-Specker coloring, and by Theorem~\ref{thm:Q} the set $\setS_3(462)$ has no Kochen-Specker coloring. Because $6 = \gcd(30,462)$, it follows from Lemma~\ref{lem:Bezout} that $\M_3(\Z[1/N])_\sym$ has no prime partial ideals. 

(2) If $N$ is even, then $\setS_d(N) \supseteq \setS_d(2)$ has no Kochen-Specker coloring by Theorem~\ref{thm:dim 4 and 6}. It follows from Lemma~\ref{lem:coloring} that $\M_d(\Z[1/N])_\sym$ has no prime partial ideals.

(3) Uncolorability of $\setS_4(2)$ as in part~(2) above implies that $\Proj(\M_4(\Z[1/2]))$ has no KS coloring. It follows from~\cite[Lemma~2.17(3)]{BMR} that $\Proj(\M_6(Z[1/2])$ also has no KS coloring. 
As shown in Theorem~\ref{thm:dim 4 and 6}, $\setS_6(3)$ also has no KS coloring. Because $1 = \gcd(2,3)$, Lemmas~\ref{lem:coloring} and~\ref{lem:Bezout} imply that $\M_6(\Z)_\sym$ has no prime partial ideals. The claim for $\M_d(\Z)_\sym$ with $d \geq 6$ now follows from~\cite[Lemma~2.17(4)]{BMR}.
\end{proof}

\separate

We are finally ready to prove our main result from Section~\ref{sec:intro}.

\begin{proof}[Proof of Theorem~\ref{thm:map to commutative}]
(1) 
First suppose that $6 \mid N$. By Theorem~\ref{thm:prime} and~\cite[Lemma~2.4]{BMR}, there is no morphism of partial rings $\M_3(\Z[1/N])_\sym \to C$ for any nonzero commutative ring $C$.

Now suppose that $p \notdivides N$ for either $p = 2,3$. Then the (unique) ring homomorphism $\Z[1/N] \to \F_p$ induces a morphism of partial rings $\M_3(\Z[1/N])_\sym \to \M_3(\F_p)_\sym$. There exists a morphism of partial rings $\M_3(\F_p)_\sym \to \F_{p^6}$ as in the proof of~\cite[Theorem~3.5]{BMR}. Composing these maps gives a morphism of partial rings
\[
\M_3(\Z[1/N])_\sym \to \F_{p^6}
\] 
to a commutative ring, proving the claim.

(2--3) The nonexistence of morphisms of partial rings in these cases follows directly from~\cite[Lemma~2.4]{BMR} and Theorem~\ref{thm:prime}.
\end{proof}

In Theorem~\ref{thm:map to commutative}(2), there is the remaining question of whether one can construct an algebraic hidden state for $\M_d(\Z[1/2])_\sym$ in dimensions $d = 4,5$. 
The construction in the case of dimension~3 began with a coloring $\Proj(\M_3(\F_3))$, which was extended using~\cite[Appendix~A]{BMR} to a morphism of partial rings to an extension of $\F_3$. But that method of extension relied crucially on the 3-dimensionality of the situation. 
It is conceivable that a similar method could leverage a coloring such as that of Proposition~\ref{prop:dim 4 and 5} to define a morphism of partial rings. But such an approach would still require some fresh insights, so we do not pursue it further in this paper.

\separate

The results of Section~\ref{sec:uncolorable} also allow us to
investigate KS colorability of projections over some rings of number-theoretic interest. 
This was explored to some extent in~\cite{BMR}, including for finite fields of prime order. We are now able to extend this treatment to $p$-adic numbers and localizations of the ring of integers in a unified way. 
For a prime $p$, let $\F_p = \Z/p\Z$ denote the field with~$p$ elements, let $\Z_p$ denote the ring of $p$-adic integers, and let $\Z_{(p)} = \{a/b \in \Q \mid a,b \in \Z, \, p \notdivides b\}$ denote the localization of the integers at the prime ideal $(p)$. 

\begin{corollary}\label{cor:projections}
Let $p$ be a prime, and let $R$ denote any of the rings $\F_p$, $\Z_{(p)}$, or $\Z_p$. 
Then:
\begin{enumerate}
\item  $\Proj(\M_3(R))$ has a Kochen-Specker coloring if and only if $p = 2,3$.
\item $\Proj(\M_4(R))$ has a Kochen-Specker coloring if and only if $p = 2$.
\item $\Proj(\M_5(R))$ has no Kochen-Specker coloring for $p > 2$.
\item For $d \geq 6$, $\Proj(\M_d(R))$ has no Kochen-Specker coloring.
\end{enumerate}
\end{corollary}

\begin{proof}
(1) For each prime $p$ we have injective and surjective ring homomorphisms $\Z_{(p)} \hookrightarrow \Z_p \twoheadrightarrow \F_p$. These induce homomorphisms between their corresponding matrix rings, which restrict to morphisms of partial Boolean algebras
\[
\Proj(\M_3(\Z_{(p)})) \to \Proj(\M_3(\Z_p)) \to \Proj(\M_3(\F_p)).
\]

For $p = 2,3$ it follows from~\cite[Theorem~3.4]{BMR} that $\Proj(\M_3(\F_p))$ has a Kochen-Specker coloring. By lifting along the morphism $\Proj(\M_3(\Z_p)) \to \Proj(\M_3(\F_p))$ (see also~\cite[Theorem~2.13]{BMR}), we obtain a coloring of the $3 \times 3$ projection matrices over $\Z_p$ and $\Z_{(p)}$.

For $p = 5$, we have $\Z[1/462] \subseteq \Z_{(5)}$ and for $p > 5$ we have $\Z[1/30] \subseteq \Z_{(p)}$. In light of Lemma~\ref{lem:coloring}, the $3 \times 3$ projection matrices over $\Z[1/462]$ are uncolorable by Theorem~\ref{thm:Q}, while those over $\Z[1/30]$ are uncolorable by~\cite{Bub:Schutte}. In either case, we can embed an uncolorable set in $\Proj(\M_3(\Z_{(p)}))$, and the image of this set remains uncolorable over $\Z_p$ and $\F_p$.

(2--4) For $p > 2$ we have $\Z[1/2] \subseteq \Z_{(p)}$, and if $d \geq 4$ then $\Proj(\M_d(\Z[1/2]))$ is uncolorable by Theorem~\ref{thm:dim 4 and 6} and Lemma~\ref{lem:coloring}. The proof of~(2) is completed by noting that $\Proj(\M_4(\F_2))$ is colorable according to Table~\ref{tab:projections}, from which $\Proj(\M_4(R))$ inherits a coloring as in part~(1). The proof of~(4) is completed by noting that $\Z[1/3] \subseteq \Z_{(2)}$ and again combining Theorem~\ref{thm:dim 4 and 6} with Lemma~\ref{lem:coloring}.
\end{proof}

\separate

In spite of the fact that Theorems~\ref{thm:map to commutative} and~\ref{thm:prime} give precise characterizations for all values of $N$ in dimensions $d = 3$ and $d \geq 6$, Question~\ref{q:colorable} remains wide open for those same dimensions. (It seems almost an accident that the question can be answered precisely for $d = 4,5$ as in Corollary~\ref{cor:dim 4 and 5}.) We close with some observations related to this number-theoretic coloring problem.

If $\setS_3(6)$ has no KS coloring, this would provide a clean, parallel answer to the question in the sense that $\setS_3(N)$ would be KS uncolorable if and only if $6 \mid N$. 
It would also provide a more straightforward proof of Theorem~\ref{thm:prime} (and thereby Theorem~\ref{thm:map to commutative}) without the need to refer to both of the sets $\setS_3(30)$ and $\setS_3(462)$.

However, our computational exploration (Table~\ref{tab:dim 3}) has failed to produce an uncolorable set in $\setS_3(6)$ to very large norm,  suggesting that 
perhaps this set has a KS coloring.
If $\setS_3(6)$ turns out to be KS colorable, it is far from clear what property would separate those values of $N$ that produce KS colorable sets from those that do not. For instance, it was shown by Salt~\cite[Table~3.2]{Salt} that $\setS_3(714)$ is also KS uncolorable, where $ 714 = 2 \cdot 3 \cdot 7 \cdot 17$.  

Below are some naive but seemingly difficult questions that highlight how much mystery remains in Question~\ref{q:colorable}:
\begin{itemize}
\item Is there any value of $d$ such that there exist infinitely many squarefree integers $\{N_i\}_{i=1}^\infty$ with $\setS_d(N_i)$ uncolorable and $N_i \nmid N_j$ for $i \neq j$?
\item Given any integer $N \geq 2$, does there exist $d$ sufficiently large such that $\setS_d(N)$ is KS uncolorable? 
\end{itemize}
A resolution of either of these questions would undoubtedly require intriguing number-theoretic insights.

\bibliographystyle{amsplain}
\bibliography{ks462-v4}
\end{document}